# A way to prove the irrationality of ζ(4)

Dirk Huylebrouck
Faculty for Architecture, KULeuven
Paleizenstraat 65-67
Brussels, 1030, BELGIUM
E-mail: Huylebrouck@gmail.com

### Abstract

*The proof of the irrationality of ζ(5) is a long standing open problem, but here only the case of ζ(4) = $\pi^4/90$ is considered. The present paper suggests an approach for the irrationality of ζ(4) along the lines of those known for proving the irrationality of ζ(3).*

## 1. Proving ζ(2), ζ(3) and ζ(4) are irrational

In 1979, F. Beukers gave an easy version of F. Apéry's proofs for the irrationality of

$$\zeta(2) = 1 + \frac{1}{2^2} + \frac{1}{3^2} \ldots \quad \text{and} \quad \zeta(3) = 1 + \frac{1}{2^3} + \frac{1}{3^3} \ldots \quad (\text{see [1]}).$$

In 1998, S. Miller modified it into a still easier proof for the irrationality of ζ(3) (see [4]). In 2001, a summary of the proofs for the irrationality of π, ln2, ζ(2) and ζ(3) was welcomed as the lack of progress in this field justified a new impulse (see [2]).

Beukers' proof of the irrationality of ζ(2) first shows that

$$0 < \left| \int_0^1 \int_0^1 \frac{(x(1-x)y(1-y))^n}{(1-xy)^{n+1}} dx dy \right| = \left| \frac{R_n + S_n \zeta(2)}{T_n} \right| \quad \text{for any } n \in \mathbb{N} \text{ and integers } R_n, S_n \text{ and } T_n.$$

Moreover, $T_n = \text{LCM}(1^2, 2^2, \ldots, n^2)$ and following a result explained in [2], $T_n \leq e^{2.01n}$ for large $n$.

Now $\left| \int_0^1 \int_0^1 \frac{(x(1-x)y(1-y))^n}{(1-xy)^{n+1}} dx dy \right| \leq \left| \int_0^1 \int_0^1 \left( \left( \frac{-1+\sqrt{5}}{2} \right)^5 \right)^n \frac{1}{(1-xy)} dx dy \right| \leq \left| \left( \frac{-1+\sqrt{5}}{2} \right)^{5n} \zeta(2) \right|.$

Thus, $0 < |R_n + S_n \zeta(2)| \leq \left| T_n \left( \frac{-1+\sqrt{5}}{2} \right)^{5n} \zeta(2) \right| \leq \left| \left( e^{2.01} \cdot \left( \frac{-1+\sqrt{5}}{2} \right)^5 \right)^n \zeta(2) \right| \to 0$ for large $n$. This is

not possible unless ζ(2) is irrational.

This proof of the irrationality of ζ(3) goes in a similar way. First, it is shown that

$$0 < \left| \int_0^1 \int_0^1 \int_0^1 \frac{(x(1-x)y(1-y)z(1-z))^n}{(1-(1-xy)z)^{n+1}} dx\,dy\,dz \right| = \left| \frac{R_n + S_n \zeta(3)}{T_n} \right| \text{ for any } n \in \mathbb{N} \text{ and integers } R_n, S_n$$

and $T_n$. Now, $T_n = \mathrm{LCM}(1^3, 2^3, \ldots, n^3)$ and following the result explained in [2], $T_n \leq e^{3.01n}$ for large $n$.

Now $\left| \int_0^1 \int_0^1 \int_0^1 \frac{(x(1-x)y(1-y)z(1-z))^n}{(1-(1-xy)z)^{n+1}} dx\,dy\,dz \right| \leq \left| \int_0^1 \int_0^1 \int_0^1 \left( \left(-1+\sqrt{2}\right)^4 \right)^n \frac{1}{(1-(1-xy)z)} dx\,dy\,dz \right|$

$$\leq \left| \left( \left(-1+\sqrt{2}\right)^4 \right)^n \zeta(3) \right|.$$

Thus, $0 < |R_n + S_n \zeta(3)| \leq \left| T_n \left( \left(-1+\sqrt{2}\right)^4 \right)^n \zeta(3) \right| \leq \left| \left( e^{3.01} \cdot \left(-1+\sqrt{2}\right)^4 \right)^n \zeta(3) \right| \to 0$ for large $n$. This is not possible unless $\zeta(3)$ is irrational.

We now suggest that a proof of the irrationality of $\zeta(4)$ could go as follows. First, it should be shown that

$$0 < \left| \int_0^1 \int_0^1 \int_0^1 \int_0^1 \frac{(x(1-x)y(1-y)z(1-z)w(1-w))^n (1-xy)^{2n+1}}{((1-(1-xy)z)(1-(1-xy)w))^{n+1}} dx\,dy\,dz\,dw \right| = \left| \frac{R_n + S_n \zeta(4)}{T_n} \right| \text{ for any}$$

$n \in \mathbb{N}$ and integers $R_n$, $S_n$ and $T_n$, where $T_n = \mathrm{LCM}(1^4, 2^4, \ldots, n^4)$. Again following the result explained in [2], $T_n \leq e^{4.01n}$ for large $n$.

Now $\left| \int_0^1 \int_0^1 \int_0^1 \int_0^1 \frac{(x(1-x)y(1-y)z(1-z)w(1-w))^n (1-xy)^{2n+1}}{((1-(1-xy)z)(1-(1-xy)w))^{n+1}} dx\,dy\,dz\,dw \right|$

$$\leq \left| \int_0^1 \int_0^1 \int_0^1 \int_0^1 \left( \frac{(-7+\sqrt{17})^4 (-3+\sqrt{17})^2}{256(1+\sqrt{17})^2} \right)^n \frac{1-xy}{(1-(1-xy)z)(1-(1-xy)w)} dx\,dy\,dz\,dw \right|$$

$$\leq \left| \left( \frac{(-7+\sqrt{17})^4 (-3+\sqrt{17})^2}{256(1+\sqrt{17})^2} \right)^n \zeta(4) \right|.$$

Thus, $0 < |R_n + S_n \zeta(4)| \leq \left| T_n \left( \frac{(-7+\sqrt{17})^4 (-3+\sqrt{17})^2}{256(1+\sqrt{17})^2} \right)^n \zeta(3) \right|$

$$\leq \left| \left( e^{4.01} \cdot \frac{(-7+\sqrt{17})^4 (-3+\sqrt{17})^2}{256(1+\sqrt{17})^2} \right)^n \zeta(4) \right| \to 0 \text{ for large } n.$$

This is not possible unless $\zeta(4)$ is irrational.

## 2. The missing part

The step $\left| \int_0^1 \int_0^1 \int_0^1 \int_0^1 \frac{(x(1-x)y(1-y)z(1-z)w(1-w))^n (1-xy)^{2n+1}}{((1-(1-xy)z)(1-(1-xy)w))^{n+1}} dx\,dy\,dz\,dw \right| = \left| \frac{R_n + S_n \zeta(4)}{T_n} \right|$ is missing and will only be shown for $n = 0, 1, 2$. In a previous paper (see [2]), it was shown that

$$\left| \int_0^1\int_0^1\int_0^1\int_0^1 \frac{(x(1-x)y(1-y)z(1-z)w(1-w))^n \cdot (1-xy))^{n+1}}{((1-(1-xy)z)(1-(1-xy)w))^{n+1}} dxdydzdw \right|$$

had potential for attempting a proof for ζ(4), but in the same paper it was also pointed out this option failed since the numerator of the integral is not of the form $R_n + S_n \zeta(4)$.

Another more esthetic expression seemed promising too (see [3]):

$$\zeta(m) = \frac{1}{m-1}\int_0^1\int_0^1\int_0^1\cdots\int_0^1 \frac{1}{(1-xy)(1-xyz)\ldots(1-xyz\ldots w)} dxdydz\ldots dw,$$

so that the expression

$$\left| \int_0^1\int_0^1\int_0^1\int_0^1 \frac{(x(1-x)y(1-y)z(1-z)w(1-w))^n}{((1-xy)(1-xyz)(1-xyzw))^{n+1}} dxdydzdw \right|$$

seemed to be a good start, but again the numerator of the integral is not of the form $R_n + S_n \zeta(4)$.

However, the current integrals

$$I_n = \int_0^1\int_0^1\int_0^1\int_0^1 \frac{(x(1-x)y(1-y)z(1-z)w(1-w))^n (1-xy)^{2n+1}}{((1-(1-xy)z)(1-(1-xy)w))^{n+1}} dxdydzdw$$

seem more promising, since at least for the values n = 0, 1 and 2 they are of the required form.

Indeed, since $\frac{1}{(1-xy)^t} = \sum_{i\geq 0}^{+\infty}\binom{i+t-1}{i}(xy)^i$ it follows that

$$\int_0^1\int_0^1 \frac{x^{r+\sigma} y^{s+\sigma}}{(1-xy)^t} dydx = \sum_{i\geq 0}^{+\infty}\binom{i+t-1}{i}\int_0^1\int_0^1 x^{i+r+\sigma} y^{i+s+\sigma} dydx = \sum_{m\geq 1}^{+\infty}\binom{i+t-1}{i}\frac{1}{(m+r+\sigma)(m+s+\sigma)}$$

with m = i+1.

Thus, after a derivation with respect to σ:

$$\int_0^1\int_0^1 \frac{x^{r+\sigma} y^{s+\sigma}}{(1-xy)^t} \text{Log}(xy)dydx = -\sum_{m\geq 1}^{+\infty}\binom{i+t-1}{i}\frac{(2m+r+s+2\sigma)}{(m+r+\sigma)^2(m+s+\sigma)^2}$$

$$\int_0^1\int_0^1 \frac{x^{r+\sigma} y^{s+\sigma}}{(1-xy)^t} \text{Log}^2(xy)dydx$$
$$= 2\sum_{m\geq 1}^{+\infty}\binom{i+t-1}{i}\frac{3m^2 + 3mr + 3ms + r^2 + rs + s^2 + 6m\sigma + 3r\sigma + 3s\sigma + 3\sigma^2}{(m+r+\sigma)^3(m+s+\sigma)^3}$$

Thus, if σ = 0 and t = 2: $\int_0^1\int_0^1 \frac{x^r y^s}{(1-xy)^2} dydx = \sum_{m\geq 1}^{+\infty}\frac{m}{(m+r)(m+s)}$

$$\int_0^1\int_0^1 \frac{x^r y^s}{(1-xy)^2} \text{Log}(xy)dydx = -\sum_{m\geq 1}^{+\infty}\frac{m(2m+r+s)}{(m+r)^2(m+s)^2}$$

$$\int_0^1\int_0^1 \frac{x^r y^s}{(1-xy)^2} \text{Log}^2(xy)dydx = 2\sum_{m\geq 1}^{+\infty}\frac{m \cdot (3m^2 + 3mr + 3ms + r^2 + rs + s^2)}{(m+r)^3(m+s)^3}$$

In case σ = 0 and t = 3:

$$\int_0^1\int_0^1 \frac{x^r y^s}{(1-xy)^3} dydx = \frac{1}{2}\sum_{m\geq 1}^{+\infty}\frac{(m+1)m}{(m+r)(m+s)}$$

$$\int_0^1\int_0^1 \frac{x^r y^s}{(1-xy)^3} Log(xy) dy dx = -\frac{1}{2}\sum_{m\geq 1}^{+\infty} \frac{(m+1)m(2m+r+s)}{(m+r)^2(m+s)^2}$$

$$\int_0^1\int_0^1 \frac{x^r y^s}{(1-xy)^3} Log^2(xy) dy dx = \sum_{m\geq 1}^{+\infty} \frac{(m+1)m.(3m^2+3mr+3ms+r^2+rs+s^2)}{(m+r)^3(m+s)^3}$$

In case $\sigma = 0$ and $t = 4$:

$$\int_0^1\int_0^1 \frac{x^r y^s}{(1-xy)^4} dy dx = \frac{1}{6}\sum_{m\geq 1}^{+\infty} \frac{m(m+1)(m+2)}{(m+r)(m+s)}$$

$$\int_0^1\int_0^1 \frac{x^r y^s}{(1-xy)^4} Log(xy) dy dx = -\frac{1}{6}\sum_{m\geq 1}^{+\infty} \frac{m(m+1)(m+2)}{(m+r)^2(m+s)^2}(2m+r+s)$$

$$\int_0^1\int_0^1 \frac{x^r y^s}{(1-xy)^4} Log^2(xy) dy dx = \frac{1}{3}\sum_{m\geq 1}^{+\infty} \frac{m(m+1)(m+2)}{(m+r)^3(m+s)^3}(3m^2+3mr+3ms+r^2+rs+s^2).$$

Firstly, for $n=0$, we note that

$$I_0 = \int_0^1\int_0^1\int_0^1\int_0^1 \frac{(1-xy)}{(1-(1-xy)z)(1-(1-xy)w)} dw dz dy dx$$

$$= \int_0^1\int_0^1 \frac{1}{(1-xy)} Log^2(xy) dy dx = 2\sum_{m\geq 1}^{+\infty} \binom{i+1-1}{i} \frac{3m^2}{m^3 m^3} = 6\zeta(4).$$

Secondly, for $n=1$, we have to compute

$$I_1 = \int_0^1\int_0^1\int_0^1\int_0^1 \frac{(x(1-x)y(1-y)z(1-z)w(1-w))^1 (1-xy)^3}{(1-(1-xy)z)^2(1-(1-xy)w)^2} dw dz dy dx$$

We do the integration with respect to $w$ and $z$ first, which goes with problem using a standard math software:

$$I_1 = \int_0^1\int_0^1 \frac{x(1-x)y(1-y)(2(1-xy)+(1+xy)Log(xy))^2}{(1-xy)^3} dy dx$$

$$= \int_0^1\int_0^1 \frac{x(1-x)y(1-y)(4(1-xy)^2+4(1-x^2y^2)Log(xy)+(1+xy)^2 Log^2(xy))}{(1-xy)^3} dy dx$$

This integral can be computed in three parts so that two parts can be computed straightforwardly using a standard math software:

$$I_{1a} = \int_0^1\int_0^1 \frac{x(1-x)y(1-y).4}{(1-xy)} dy dx = -13 + 8\zeta(2)$$

$$I_{1b} = \int_0^1\int_0^1 \frac{x(-1+x)y(-1+y).(-4).(1+xy)}{(-1+xy)^2} Log(xy) dy dx = -51 -16\zeta(2) - 64\zeta(3)$$

The third part,

$$I_{1c} = \int_0^1\int_0^1 \frac{x(1-x)y(1-y).(1+xy)^2}{(1-xy)^3} Log^2(xy) dy dx$$

resists the software, and thus we expand the numerator in order to compute it term by term:

$$I_{1c} = \int_0^1\int_0^1 \frac{xy - x^2y - xy^2 + 3x^2y^2 - 2x^3y^2 - 2x^2y^3 + 3x^3y^3 - x^4y^3 - x^3y^4 + x^4y^4}{(1-xy)^3} Log^2(xy)\,dy\,dx$$

Using the above expressions:

$$I_{1c} = \sum_{m\geq 1}^{+\infty}\left(\frac{-3}{(m+1)^3} + \frac{5}{(m+1)^2} + \frac{-2}{m+1} + \frac{18}{(m+2)^4} + \frac{-31}{(m+2)^3} + \frac{15}{(m+2)^2} + \frac{-2}{m+2}\right.$$

$$\left. + \frac{54}{(m+3)^4} + \frac{-33}{(m+3)^3} + \frac{-1}{(m+3)^2} + \frac{2}{m+3} + \frac{36}{(m+4)^4} + \frac{3}{(m+4)^3} + \frac{-11}{(m+4)^2} + \frac{2}{m+4}\right)$$

$$= \sum_{n\geq 1}^{+\infty}\left(\frac{-3}{n^3} + \frac{5}{n^2} + \frac{-2}{n} + \frac{18}{n^4} + \frac{-31}{n^3} + \frac{15}{n^2} + \frac{-2}{n} + \frac{54}{n^4} + \frac{-33}{n^3} + \frac{-1}{n^2} + \frac{2}{n} + \frac{36}{n^4} + \frac{3}{n^3} + \frac{11}{n^2} + \frac{2}{n}\right)$$

$$- 3(-1) + 5(-1) - 2(-1) + 18\left(-\frac{1}{2^4} - 1\right) - 31\left(-\frac{1}{2^3} - 1\right) + 15\left(-\frac{1}{2^2} - 1\right) - 2\left(-\frac{1}{2} - 1\right)$$

$$+ 54\left(-\frac{1}{3^4} - \frac{1}{2^4} - 1\right) - 33\left(-\frac{1}{3^3} - \frac{1}{2^3} - 1\right) - 1\left(-\frac{1}{3^2} - \frac{1}{2^2} - 1\right) + 2\left(-\frac{1}{3} - \frac{1}{2} - 1\right) +$$

$$36\left(-\frac{1}{4^4} - \frac{1}{3^4} - \frac{1}{2^4} - 1\right) + 3\left(-\frac{1}{4^3} - \frac{1}{3^3} - \frac{1}{2^3} - 1\right) - 11\left(-\frac{1}{4^2} - \frac{1}{3^2} - \frac{1}{2^2} - 1\right) + 2\left(-\frac{1}{4} - \frac{1}{3} - \frac{1}{2} - 1\right)$$

$$= 108\zeta(4) + 64\zeta(3) + 8\zeta(2) - 423/8$$

And thus

$I_1 = I_{1a} + I_{1b} + I_{1c} = (-13 + 8\zeta(2)) + (-51 - 16\zeta(2) - 64\zeta(3)) + (72\zeta(4) + 64\zeta(3) + 8\zeta(2) - 423/8) = 108\zeta(4) - 935/8$.

Note that $8*108\zeta(4) - 935 = 0.127274... < 1$

Thirdly, for computing $I_2$ we again do the integration with respect to $w$ and $z$ first:

$$I_2 = \int_0^1\int_0^1\int_0^1\int_0^1 \frac{(x(1-x)y(1-y)z(1-z)w(1-w))^2(1-xy)^5}{(1-(1-xy)z)^2(1-(1-xy)w)^3}\,dw\,dz\,dy\,dx$$

$$= \int_0^1\int_0^1 \frac{x^2(1-x)^2 y^2(1-y)^2(3(1-x^2y^2) + (1+4xy+x^2y^2)Log(xy))^2}{(1-xy)^5}\,dy\,dx$$

This integral can be cut in three parts and the first can be computed straightforwardly using standard software:

$$I_{2a} = \int_0^1\int_0^1 \frac{x^2(1-x)^2 y^2(1-y)^2(3(1-x^2y^2))^2}{(1-xy)^5}\,dy\,dx = \frac{21}{16}(-1737 + 1056\zeta(2))$$

The second and third parts resisted our software, and so they are computed using series

$$I_{2b} = \int_0^1\int_0^1 \frac{x^2(1-x)^2 y^2(1-y)^2 6(1+xy)(1+4xy+x^2y^2)Log(xy)}{(1-xy)^4}\,dy\,dx$$

$$= 6\int_0^1\int_0^1 \frac{x^2y^2 - 2x^3y^2 + x^4y^2 - ... + x^7y^7}{(1-xy)^4} Log(xy)\,dy\,dx$$

and

$$I_{2c} = \int_0^1\int_0^1 \frac{x^2y^2 - 2x^3y^2 + x^4y^2 - ... + x^8y^8}{(1-xy)^5} Log(xy)^2\,dy\,dx$$

Using the above expressions:

$$= \sum_{m\geq 1}^{+\infty}\frac{1}{2(m+2)^3} - \sum_{m\geq 1}^{+\infty}\frac{19}{4(m+2)^2} + \sum_{m\geq 1}^{+\infty}\frac{47}{4(m+2)} + \ldots$$

$$\ldots + \sum_{m\geq 1}^{+\infty}\frac{420}{(m+8)^4} + \sum_{m\geq 1}^{+\infty}\frac{307}{2(m+8)^3} + \sum_{m\geq 1}^{+\infty}\frac{815}{4(m+8)^2} + \sum_{m\geq 1}^{+\infty}\frac{225}{4(m+8)}$$

$$= \frac{1}{2}\left(\sum_{n\geq 1}^{+\infty}\frac{1}{n^3} - \frac{1}{2^3} - 1\right) - \frac{19}{4}\left(\sum_{n\geq 1}^{+\infty}\frac{1}{n^2} - \frac{1}{2^2} - 1\right) + \frac{47}{4}\left(\sum_{n\geq 1}^{+\infty}\frac{1}{n} - 2 - 1\right) + \ldots$$

$$\ldots + \frac{225}{4}\left(\sum_{n\geq 1}^{+\infty}\frac{1}{n} - \frac{1}{8} - \frac{1}{7} - \frac{1}{6} - \frac{1}{5} - \frac{1}{4} - \frac{1}{3} - \frac{1}{2} - 1\right)$$

$$= -\frac{31306541}{3456} + 10476\sum_{n\geq 1}^{+\infty}\frac{1}{n^4} - 1386\sum_{n\geq 1}^{+\infty}\frac{1}{n^2}.$$

So that $I_2 = I_{2a} + (I_{2b} + I_{2c}) = \frac{21}{16}(-1737 + 1056\zeta(2)) - \frac{31306541}{3456} + 10476\zeta(4) - 1386\zeta(2)$

$$= -\frac{39185573}{3456} + 10476\zeta(4)$$

Again, we note that $3456*10476\zeta(4) - 39185573 = 3456*0.000082932 = 0.286613 < 1$

Thus, the integrals

$$I_n = \int_0^1\int_0^1\int_0^1\int_0^1 \frac{(x(1-x)y(1-y)z(1-z)w(1-w))^n(1-xy)^{2n+1}}{((1-(1-xy)z)(1-(1-xy)w))^{n+1}}dxdydzdw$$

seem very promising indeed.